# Fibonacci polynomials, generalized Stirling numbers, and Bernoulli, Genocchi and tangent numbers

Johann Cigler

johann.cigler@univie.ac.at

**Abstract**

We study matrices which transform the sequence of Fibonacci or Lucas polynomials with even index to those with odd index and vice versa. They turn out to be intimately related to generalized Stirling numbers and to Bernoulli, Genocchi and tangent numbers and give rise to various identities between these numbers. There is also a close connection with the Akiyama-Tanigawa algorithm. Since such numbers have been extensively studied it is possible that some of these results are already hidden in the literature. I would be very grateful for such information.

**0. Introduction**

Let $(F_n(s))$ be the sequence of Fibonacci polynomials. The generating functions of the subsequences of polynomials with even and odd index are related by

$$\sum_{n\geq 0} \frac{F_{2n+2}(s)}{(2n+2)!} z^{2n+2} = \frac{e^z - 1}{e^z + 1} \sum_{n\geq 0} \frac{F_{2n+1}(s)}{(2n+1)!} z^{2n+1}.$$

If we associate with each sequence $x = (x_n)_{n\geq 0}$ the sequence $y = (y_n)_{n\geq 0}$ defined by

$$\sum_{n\geq 0} \frac{y_n}{(2n+2)!} z^{2n+2} = \frac{e^z - 1}{e^z + 1} \sum_{n\geq 0} \frac{x_n}{(2n+1)!} z^{2n+1}$$

we obtain an infinite triangular matrix $A = (a(n,k))$ such that

$$y_n = \sum_{k=0}^{n} a(n,k) x_k.$$

The matrix $A$ begins with

$$\begin{pmatrix} 1 & 0 & 0 & 0 & 0 \\ -1 & 2 & 0 & 0 & 0 \\ 3 & -5 & 3 & 0 & 0 \\ -17 & 28 & -14 & 4 & 0 \\ 155 & -255 & 126 & -30 & 5 \end{pmatrix}.$$

The above mentioned identity of generating functions implies that

$$A = \left([j \leq i]\binom{2i+1-j}{j}\right)\left([j \leq i]\binom{2i-j}{j}\right)^{-1}$$

and that

$$a(n,k) = (-1)^{n-k} \frac{1}{2k+1}\binom{2n+2}{2k} G_{2n-2k+2}$$

if we denote by $G_{2n}$ the positive Genocchi numbers.

The starting point of this paper was the observation that the eigenvectors of $A$ are sequences of central factorial numbers $T(n,k)$. More precisely the eigenvector corresponding to the eigenvalue $k+1$, $k \in \mathbb{N}$, is $(T(n+1,k+1))_{n\geq 0}$. In other words $A$ can be written in the form

$$A = (T(i+1,j+1)) D (T(i+1,j+1))^{-1},$$

where $D$ denotes the diagonal matrix with entries $1, 2, 3, \cdots$. This fact is equivalent with the identity

$$\frac{e^z - 1}{e^z + 1} \sum_{n\geq 0} \frac{T(n+1,k+1)}{(2n+1)!} z^{2n+1} = (k+1) \sum_{n\geq 0} \frac{T(n+1,k+1)}{(2n+2)!} z^{2n+2}.$$

If we define linear functionals $\varphi_{k+1}$ on the vector space of polynomials in $s$ by

$$\varphi_{k+1}(F_{2n+1}(s)) = T(n+1,k+1)$$

then

$$\varphi_{k+1}(F_{2n+2}(s)) = (k+1)T(n+1,k+1)$$

and

$$\varphi_{k+1}(s^n) = LS(n,k),$$

where $LS(n,k)$ are the Legendre-Stirling numbers which have been studied by G.E. Andrews and Littlejohn [2] and by Y. Gelineau and J. Zeng [8].

These facts lead to some interesting identities such as

$$\sum_{k=1}^{n} (-1)^{n-k} k \left((k-1)!\right)^2 T(n,k) = G_{2n}$$

or

$$\sum_{k=0}^{n}(-1)^k \frac{(k!)^2}{k+1}T(n+1,k+1) = (2n+1)B_{2n}.$$

To all these formulae the Akiyama-Tanigawa algorithm is applicable.

Another property of $A$ is the fact that it can be written as

$$A = \left([j \leq i]\binom{i+1}{2i-2j}\right)^{-1}\left([j \leq i]\binom{i+1}{2i-2j+1}\right).$$

This leads among other things to new proofs of Seidel's identity for Genocchi numbers and of Kaneko's identity.

The aim of the present paper is to give a systematic account of these and related results.

**1. Some well-known facts**

We consider (a variant of) the Fibonacci polynomials defined by

$$F_n(s) = \sum_{k=0}^{\lfloor\frac{n-1}{2}\rfloor}\binom{n-1-k}{k}s^k. \tag{1.1}$$

They satisfy the recursion $F_n(s) = F_{n-1}(s) + sF_{n-2}(s)$ with initial values $F_0(s) = 0$ and $F_1(s) = 1$.

The first terms are

$(0, 1, 1, 1+s, 1+2s, 1+3s+s^2, 1+4s+3s^2, 1+5s+6s^2+s^3, 1+6s+10s^2+4s^3, \dots)$.

The corresponding Lucas polynomials are defined by

$$L_n(s) = \sum_{k=0}^{\lfloor\frac{n}{2}\rfloor}\frac{n}{n-k}\binom{n-k}{k}s^k \tag{1.2}$$

and satisfy the same recurrence as the Fibonacci polynomials, but with initial values $L_0(s) = 2$ and $L_1(s) = 1$. The first terms of this sequence are

$(2, 1, 1+2s, 1+3s, 1+4s+2s^2, 1+5s+5s^2, 1+6s+9s^2+2s^3, \dots)$.

Let

$$\alpha = \frac{1+\sqrt{1+4s}}{2} \tag{1.3}$$

and

$$\beta = \frac{1-\sqrt{1+4s}}{2} \tag{1.4}$$

be the roots of the equation $z^2 - z - s = 0$.

Then for $s \neq -\frac{1}{4}$ the well-known Binet formulae give

$$F_n(s) = \frac{\alpha^n - \beta^n}{\alpha - \beta} \tag{1.5}$$

and

$$L_n(s) = \alpha^n + \beta^n. \tag{1.6}$$

For $s = -\frac{1}{4}$ we have instead

$$F_n\left(-\frac{1}{4}\right) = \frac{n}{2^{n-1}}. \tag{1.7}$$

It is also well known that $L_n(s) = F_{n+1}(s) + sF_{n-1}(s)$. This follows immediately from $s = -\alpha\beta$ and $\alpha^{n+1} - \beta^{n+1} - \alpha\beta(\alpha^{n-1} - \beta^{n-1}) = \alpha^n(\alpha - \beta) + \beta^n(\alpha - \beta)$.

Since $\deg F_{2n+1}(s) = \deg F_{2n+2}(s) = \deg L_{2n}(s) = \deg L_{2n+1}(s) = n$ each of the sets $\{F_{2n+1}(s)\}, \{F_{2n+2}(s)\}, \{L_{2n}(s)\}, \{L_{2n+1}(s)\}$ is a basis for the vector space of polynomials in $s$.

The Genocchi numbers $(G_{2n})_{n \geq 1} = (1,1,3,17,155,2073,38227,929569,\cdots)$ (cf. OEIS A110501) and their relatives $(g_n)_{n \geq 1} = (1,-1,0,1,0,-3,0,17,0,-155,\cdots)$ (cf. OEIS A036968) are defined by their exponential generating function

$$\frac{2z}{1+e^z} = z + z\frac{1-e^z}{1+e^z} = \sum_{n \geq 1} g_n \frac{z^n}{n!} = z + \sum_{n \geq 1} (-1)^n G_{2n} \frac{z^{2n}}{(2n)!}, \tag{1.8}$$

and the tangent numbers $(T_{2n+1}) = (1, 2, 16, 272, 7936, 353792, 22368256, \cdots)$ (cf. OEIS A000182) are defined by

$$\frac{e^{2z} - 1}{e^{2z} + 1} = \sum_{k \geq 0} (-1)^k T_{2k+1} \frac{z^{2k+1}}{(2k+1)!}. \tag{1.9}$$

Note that by comparing coefficients we get

$$(-1)^{k-1} \frac{g_{2k+2}}{2k+2} = \frac{G_{2k+2}}{2k+2} = \frac{T_{2k+1}}{2^{2k+1}}. \tag{1.10}$$

It is also well-known that $G_{2n} = (-1)^n 2(1 - 2^{2n}) B_{2n}$, where $(B_n)$ is the sequence of Bernoulli numbers defined by

$$B_n = \sum_{k=0}^{n} \binom{n}{k} B_k \tag{1.11}$$

for $n > 1$ with $B_0 = 1$. The list of Bernoulli numbers begins with

$$1, -\frac{1}{2}, \frac{1}{6}, 0, -\frac{1}{30}, 0, \frac{1}{42}, 0, -\frac{1}{30}, 0, \frac{5}{66}, 0, -\frac{691}{2730}, 0, \frac{7}{6}, 0, -\frac{3617}{510}, \ldots$$

For later uses let us also recall the generating function of the Bernoulli numbers

$$\sum_{n \geq 0} B_n \frac{z^n}{n!} = \frac{z}{e^z - 1}. \tag{1.12}$$

## 2. Connection constants

The following theorem gives an explicit computation of some basis transformations.

**Theorem 2.1**

*The bases $(F_{2n+2})$ and $(F_{2n+1})$ are connected by*

$$F_{2n+2}(s) = \sum_{k=0}^{n} (-1)^{n-k} \frac{G_{2n-2k+2}}{2k+1} \binom{2n+2}{2k} F_{2k+1}, \tag{2.1}$$

*and*

$$F_{2n+1}(s) = \sum_{k=0}^{n} \binom{2n+1}{2k+1} \frac{B_{2n-2k}}{k+1} F_{2k+2}(s). \tag{2.2}$$

*The bases consisting of Lucas polynomials are connected by*

$$L_{2n+1}(s) = \sum_{k=0}^{n}(-1)^{n-k}\frac{T_{2n-2k+1}}{2^{2n-2k+1}}\binom{2n+1}{2k}L_{2k}(s) = \sum_{k=0}^{n}(-1)^{n-k}\frac{G_{2n-2k+2}}{2n-2k+2}\binom{2n+1}{2k}L_{2k}(s) \quad (2.3)$$

*and*

$$L_{2n}(s) = 2\sum_{j=0}^{n}\binom{2n}{2j}\frac{B_{2n-2j}}{2j+1}L_{2j+1}(s). \quad (2.4)$$

**Proof**

1) Since $\alpha + \beta = 1$ we have

$$e^{\alpha z} - e^{\beta z} = e^{(1-\beta)z} - e^{(1-\alpha)z} = -e^{z}\left(e^{-\alpha z} - e^{-\beta z}\right). \quad (2.5)$$

This gives

$$\sum_{n\geq 0}\frac{F_n(s)}{n!}z^n = -e^{z}\sum_{n\geq 0}\frac{F_n(s)}{n!}(-z)^n \quad (2.6)$$

or

$$\frac{1+e^{-z}}{2}\sum_{n}F_n(s)\frac{z^n}{n!} = \sum_{n}F_{2n+1}(s)\frac{z^{2n+1}}{(2n+1)!} \quad (2.7)$$

This reduces our task to finding the matrix which corresponds to the operator "multiplication by $\frac{1+e^{-z}}{2}$" of the corresponding exponential generating functions.

Comparing coefficients we see that for $x = (x_n)_{n\geq 0}$ and $y = (y_n)_{n\geq 0}$

$$\frac{1+e^{-z}}{2}\sum_{n}x_n\frac{z^n}{n!} = \sum_{n}y_n\frac{z^n}{n!} \quad (2.8)$$

is equivalent with $y = Mx$

where $M = (m(i,j))$ with $m(n,k) = (-1)^{n-k}\frac{1}{2}\binom{n}{k}$ for $k < n$, $m(n,n) = 1$ and $m(n,k) = 0$ for $k > n$.

The inverse of (2.8) is

$$\sum_n x_n \frac{z^n}{n!} = \frac{2}{1+e^{-z}} \sum_n y_n \frac{z^n}{n!} = \sum_j g_{j+1} \frac{(-z)^j}{(j+1)j!} \sum_k y_k \frac{z^k}{k!} = \sum_n \frac{z^n}{n!} \sum_{k=0}^n (-1)^{n-k} \frac{g_{n-k+1}}{n-k+1} \binom{n}{k} y_k, \text{ thus}$$

$$x_n = \sum_{k=0}^n (-1)^{n-k} \frac{g_{n-k+1}}{n-k+1} \binom{n}{k} y_k. \tag{2.9}$$

Thus (2.7) implies

$$F_n(s) = \sum (-1)^{n-1} \frac{g_{n+2-2k}}{n+2-2k} \binom{n}{2k-1} F_{2k-1}$$

and we get as special case

$$F_{2n+2}(s) = -\sum_{k \geq 0} \frac{g_{2n-2k+2}}{2n-2k+2} \binom{2n+2}{2k+1} F_{2k+1} = \sum_{k=0}^n (-1)^{n-k} \frac{G_{2n-2k+2}}{2k+1} \binom{2n+2}{2k} F_{2k+1},$$

i.e. (2.1).

From

$$e^{\alpha z} + e^{\beta z} = e^{(1-\beta)z} + e^{(1-\alpha)z} = e^z \left( e^{-\alpha z} + e^{-\beta z} \right) \tag{2.10}$$

we get in the same way as before

$$\frac{1+e^{-z}}{2} \sum_n L_n(s) \frac{z^n}{n!} = \sum_n L_{2n}(s) \frac{z^{2n}}{(2n)!}. \tag{2.11}$$

In this case (2.9) implies (2.3).

2) Another consequence of (2.5) is

$$(1-e^{-z}) \sum_{n \geq 0} \frac{F_n(s)}{n!} z^n = 2 \sum_{n \geq 0} \frac{F_{2n}(s)}{(2n)!} z^{2n}. \tag{2.12}$$

Now

$$(1-e^{-z}) \sum_k x_k \frac{z^k}{k!} = \sum_n \frac{z^n}{n!} \sum_{j=0}^{n-1} (-1)^{n-j-1} \binom{n}{j} x_j = \sum_{n \geq 1} y_{n-1} \frac{z^n}{n!} \text{ is equivalent with}$$

$$y_n = \sum_{j=0}^{n} (-1)^{n-j} \binom{n+1}{j} x_j. \tag{2.13}$$

The inverse is

$$\sum_k x_k \frac{z^k}{k!} = \frac{1}{(1-e^{-z})} \sum_n y_{n-1} \frac{z^n}{n!} = \sum_n \frac{b_k z^{k-1}}{k!} \sum_n y_{j-1} \frac{z^j}{j!}$$

where $b(n) = B_n$ for $n \neq 1$ and $b(1) = \frac{1}{2}$. This follows from

$$\sum_n b_n \frac{z^n}{n!} = \sum_n B_n \frac{z^n}{n!} + z = \frac{z}{e^z - 1} + z = \frac{ze^z}{e^z - 1} = \frac{z}{1 - e^{-z}}.$$

This implies

$$x_n = \sum_{j=0}^{n} \binom{n}{j} \frac{b_{n-j}}{j+1} y_j. \tag{2.14}$$

From (2.12) we see that with $x_n = F_n(s)$ and $y_{2n-1} = 2F_{2n}(s), y_{2n} = 0,$ we get

$$F_{2n+1}(s) = \sum_{j=0}^{2n+1} \binom{2n+1}{j} \frac{b_{2n+1-j}}{j+1} y(j) = \sum_{j=0}^{n} \binom{2n+1}{2j+1} \frac{b_{2n-2j}}{2j+2} y(2j+1) = \sum_{j=0}^{n} \binom{2n+1}{2j+1} \frac{b_{2n-2j}}{j+1} F_{2j+2},$$

i.e. (2.2).

(2.10) implies $(1-e^{-z}) \sum_{n \geq 0} \frac{L_n(s)}{n!} z^n = 2 \sum_{n \geq 0} \frac{L_{2n+1}(s)}{(2n+1)!} z^{2n+1}.$

Therefore we get from (2.12) with $x_n = L_n(s)$ and $y_{2n} = 2L_{2n+1}(s), y_{2n+1} = 0,$

$$L_{2n}(s) = \sum_{j=0}^{n} \binom{2n}{2j} \frac{b_{2n-2j}}{2j+1} y(2j) = \sum_{j=0}^{n} \binom{2n}{2j} \frac{B_{2n-2j}}{2j+1} 2L_{2j+1}(s),$$

i.e. (2.4).

The matrix corresponding to "multiplication by $(1-e^{-z})$" of the generating functions is

$$C = \left( (-1)^{i-j} [j \leq i] \binom{i+1}{j} \right)_{i,j=0}^{\infty}. \tag{2.15}$$

Its inverse is given by

$$C^{-1} = \left( [j \leq i] \binom{i}{j} \frac{b(i-j)}{j+1} \right)_{i,j=0}^{\infty}. \tag{2.16}$$

Here $[P] = 1$ of property $P$ is true and $[P] = 0$ else.

If we apply $C^{-1}$ to the column vector $y$ with entries $y_{2n+1} = 2F_{n+2}(s)$, $y_{2n} = 0$, then the entries $x_n$ of $x = C^{-1}y$ are $x_n = F_n(s)$.

For even $n$ this is trivial because in $x_{2n} = \sum_{j=0}^{n} \binom{2n}{2j+1} \frac{b_{2n-2j-1}}{2j+2} y_{2j+1} = \binom{2n}{2n-1} \frac{F_{2n}}{2n} = F_{2n}$ the only non-vanishing term occurs for $j = n-1$. For odd $n$ we get the identities (2.2).

If we apply $C^{-1}$ to the column vector $y$ with entries $y_{2n} = 2L_n(s)$, $y_{2n+1} = 0$, then the entries $x_n$ of $x = C^{-1}y$ are $x_n = L_n(s)$. For even $n$ we get (2.4). For odd $n$ the corresponding identity is again trivial.

The matrices (2.15) and (2.16) are intimately connected with Stirling numbers.

## 3. Stirling numbers and central factorial numbers

**3.1.** In order to make the paper self-contained I shall first state some well-known facts about Stirling numbers and some generalizations (cf. [4]).

Let $w = (w(n))_{n \geq 0}$ be an increasing sequence of positive numbers. Define the $w$-Stirling numbers of the second kind $S^w(n,k)$ by

$$S^w(n,k) = S^w(n-1,k-1) + w(k)S^w(n-1,k) \tag{3.1}$$

with initial values $S^w(0,k) = [k = 0]$ and $S^w(n,0) = w(0)^n$ and the $w$-Stirling numbers of the first kind $s^w(n,k)$ by

$$s^w(n,k) = s^w(n-1,k-1) - w(n-1)s^w(n-1,k) \tag{3.2}$$

with $s^w(0,k) = [k = 0]$ and $s^w(n,0) = (-1)^n \prod_{j=0}^{n-1} w(j)$.

This is equivalent with

$$\sum_{k=0}^{n} s^w(n,k) x^k = (x - w(0))(x - w(1)) \cdots (x - w(n-1)), \tag{3.3}$$

$$\sum_{k=0}^{n} S^w(n,k) \prod_{j=0}^{k-1} (x - w(j)) = x^n \tag{3.4}$$

and

$$\sum_{n \geq k} S^w(n,k) x^n = \frac{x^k}{(1 - w(0)x)(1 - w(1)x) \cdots (1 - w(k)x)}. \tag{3.5}$$

Let

$$S^w(n) = \left(S^w(i,j)\right)_{i,j \geq 0} \tag{3.6}$$

and

$$s^w(n) = \left(s^w(i,j)\right)_{i,j \geq 0}. \tag{3.7}$$

Then it is clear that

$$s^w(n) = \left(S^w(n)\right)^{-1}. \tag{3.8}$$

The same holds for the infinite matrices $S^w = \left(S^w(i,j)\right)_{i,j \geq 0}$ and $s^w = \left(s^w(i,j)\right)_{i,j \geq 0}$.

**3.2.** For $w(n) = n$ we get the Stirling numbers $s(n,k)$ of the first kind and the Stirling numbers $S(n,k)$ of the second kind.

Note that the matrices $\left(S(i,j)\right)_{i,j \geq 1}$ and $\left(s(i,j)\right)_{i,j \geq 1}$ are the same as $\left(S^w(i,j)\right)_{i,j \geq 0}$ and $\left(s^w(i,j)\right)_{i,j \geq 0}$ respectively for $w(n) = n+1$.

Consider the infinite matrix $C = \left(c(i,j)\right)_{i,j \geq 0}$ with $c(i,j) = (-1)^{i-j} \binom{i+1}{j}$ for $j \leq i$ and $c(i,j) = 0$ for $j > i$, which has been defined in (2.15).

**Theorem 3.1**

*The matrix C can be factored in the following way:*

$$C = \left( \left( [j \leq i] \binom{i+1}{j} \right)_{i,j \geq 0} \right)^{-1} \left( [j \leq i] \binom{i}{j} \right)_{i,j \geq 0} \quad (3.9)$$

$$= \left( S(i+1, j+1) \right)_{i,j \geq 0} \left( (i+1)([i = j]) \right)_{i,j \geq 0} \left( s(i+1, j+1) \right)_{i,j \geq 0}.$$

**Proof**

This theorem can be derived from the following simple identities.

The identity

$$C \left( [j \leq i] \binom{i}{j} \right)_{i,j \geq 0} = \left( [j \leq i] \binom{i+1}{j} \right)_{i,j \geq 0} \quad (3.10)$$

is another formulation of the trivial result

$$\sum_{j=0}^{n} (-1)^{n-j} \binom{n+1}{j} \binom{j}{k} = \binom{n+1}{k} + \sum_{j=0}^{n+1} (-1)^{n-j} \binom{n+1}{j} \binom{j}{k} = \binom{n+1}{k}$$

for $k \leq n$.

The identity

$$\left( [j \leq i] \binom{i}{j} \right)_{i,j \geq 0} \left( S(i, j) \right)_{i,j \geq 0} = \left( S(i+1, j+1) \right)_{i,j \geq 0} \quad (3.11)$$

is a well-known property of the Stirling numbers and easily proved by induction:

$$\sum_{j=0}^{n} \binom{n}{j} S(j,k) = \sum_{j=0}^{n-1} \binom{n-1}{j} S(j,k) + \sum_{j=0}^{n-1} \binom{n-1}{j-1} S(j,k) = S(n,k+1) + \sum_{j=0}^{n-1} \binom{n-1}{j} S(j+1,k)$$

$$= S(n,k+1) + \sum_{j=0}^{n-1} \binom{n-1}{j} S(j,k-1) + k \sum_{j=0}^{n-1} \binom{n-1}{j} S(j,k)$$

$$= S(n,k+1) + S(n,k) + kS(n,k+1) = S(n+1,k+1).$$

Finally we have

$$\left([j\leq i]\binom{i+1}{j}\right)_{i,j\geq 0}(S(i,j))_{i,j\geq 0} = (S(i+1,j+1)(j+1))_{i,j\geq 0}. \tag{3.12}$$

This also is easily proved by induction:

$$\sum_{j=0}^{n}\binom{n+1}{j}S(j,k) = \sum_{j=0}^{n}\binom{n}{j}S(j,k) + \sum_{j=0}^{n}\binom{n}{j-1}S(j,k) = S(n+1,k+1) + \sum_{j=0}^{n-1}\binom{n}{j}S(j+1,k)$$

$$= S(n+1,k+1) + \sum_{j=0}^{n-1}\binom{n}{j}S(j,k-1) + k\sum_{j=0}^{n-1}\binom{n}{j}S(j,k)$$

$$= S(n+1,k+1) + \left(S(n+1,k) - S(n,k-1)\right) + k\left(S(n+1,k+1) - S(n,k)\right)$$

$$= (k+1)S(n+1,k+1) + \left(S(n+1,k) - S(n,k-1) - kS(n,k)\right) = (k+1)S(n+1,k+1).$$

To derive (3.9) observe that (3.10) gives

$$C = \left([j\leq i]\binom{i+1}{j}\right)_{i,j\geq 0}\left(\left([j\leq i]\binom{i}{j}\right)_{i,j\geq 0}\right)^{-1},$$

(3.12) gives

$$\left([j\leq i]\binom{i+1}{j}\right)_{i,j\geq 0} = (S(i+1,j+1)(j+1))_{i,j\geq 0}\left((S(i,j))_{i,j\geq 0}\right)^{-1}$$

and (3.11) gives

$$\left(\left([j\leq i]\binom{i}{j}\right)_{i,j\geq 0}\right)^{-1} = (S(i,j))_{i,j\geq 0}\left((S(i+1,j+1))_{i,j\geq 0}\right)^{-1}.$$

Combining these identities gives (3.9).

Another consequence of (3.11) and (3.12) is

$$\left(\left([j\leq i]\binom{i}{j}\right)_{i,j\geq 0}\right)^{-1}\left([j\leq i]\binom{i+1}{j}\right)_{i,j\geq 0} = (S(i,j))_{i,j\geq 0}((i+1)([i=j]))_{i,j\geq 0}(s(i,j))_{i,j\geq 0}. \tag{3.13}$$

**3.3.** If we choose $w(n) = n^2$ we get the central factorial numbers $t(n,k) = s^w(n,k)$ of the first kind and the central factorial numbers $T(n,k) = S^w(n,k)$ of the second kind respectively.

These numbers have been introduced in [12] with a different notation. Further results can be found in [14], Exercise 5.8.

The following tables show the upper part of the matrices of central factorial numbers. See also OEIS A036969.

$$(T(i,j))_{i,j=0}^{6} = \begin{pmatrix} 1 & 0 & 0 & 0 & 0 & 0 & 0 \\ 0 & 1 & 0 & 0 & 0 & 0 & 0 \\ 0 & 1 & 1 & 0 & 0 & 0 & 0 \\ 0 & 1 & 5 & 1 & 0 & 0 & 0 \\ 0 & 1 & 21 & 14 & 1 & 0 & 0 \\ 0 & 1 & 85 & 147 & 30 & 1 & 0 \\ 0 & 1 & 341 & 1408 & 627 & 55 & 1 \end{pmatrix}$$

$$(t(i,j))_{i,j=0}^{6} = \begin{pmatrix} 1 & 0 & 0 & 0 & 0 & 0 & 0 \\ 0 & 1 & 0 & 0 & 0 & 0 & 0 \\ 0 & -1 & 1 & 0 & 0 & 0 & 0 \\ 0 & 4 & -5 & 1 & 0 & 0 & 0 \\ 0 & -36 & 49 & -14 & 1 & 0 & 0 \\ 0 & 576 & -820 & 273 & -30 & 1 & 0 \\ 0 & -14400 & 21076 & -7645 & 1023 & -55 & 1 \end{pmatrix}.$$

**3.4** For $w(n) = n(n+1)$ we get the Legendre-Stirling numbers $LS(n,k)$ and $ls(n,k)$ studied in [2] and [8] (cf. OEIS A071951 and A129467).

$$(LS(i,j))_{i,j=0}^{6} = \begin{pmatrix} 1 & 0 & 0 & 0 & 0 & 0 & 0 \\ 0 & 1 & 0 & 0 & 0 & 0 & 0 \\ 0 & 2 & 1 & 0 & 0 & 0 & 0 \\ 0 & 4 & 8 & 1 & 0 & 0 & 0 \\ 0 & 8 & 52 & 20 & 1 & 0 & 0 \\ 0 & 16 & 320 & 292 & 40 & 1 & 0 \\ 0 & 32 & 1936 & 3824 & 1092 & 70 & 1 \end{pmatrix}$$

$$(ls(i,j))_{i,j=0}^{6} = \begin{pmatrix} 1 & 0 & 0 & 0 & 0 & 0 & 0 \\ 0 & 1 & 0 & 0 & 0 & 0 & 0 \\ 0 & -2 & 1 & 0 & 0 & 0 & 0 \\ 0 & 12 & -8 & 1 & 0 & 0 & 0 \\ 0 & -144 & 108 & -20 & 1 & 0 & 0 \\ 0 & 2880 & -2304 & 508 & -40 & 1 & 0 \\ 0 & -86400 & 72000 & -17544 & 1708 & -70 & 1 \end{pmatrix}.$$

There are some interesting relations between central factorial numbers and Legendre-Stirling numbers which are analogous to the corresponding results about Stirling numbers.

**Theorem 3.2**

$$\sum_{j=0}^{n} \binom{2n-j}{j} LS(j,k) = T(n+1,k+1) \qquad (3.14)$$

*and*

$$\sum_{j=0}^{n} \binom{2n+1-j}{j} LS(j,k) = (k+1)T(n+1,k+1) \qquad (3.15)$$

**Proof**

For $n < k+1$ all terms vanish. For $n = k+1$ we have

$$\sum_{j=0}^{k} \binom{2k+1-j}{j} LS(j,k) = \binom{2k+1-k}{k} = (k+1) = (k+1)T(k+1,k+1)$$

and

$$\sum_{j=0}^{k} \binom{2k-j}{j} LS(j,k) = \binom{k}{k} = 1 = T(k+1,k+1).$$

Now assume that (3.14) and (3.15) are already known up to $n-1$. Then

$$\sum_{j=0}^{n} \binom{2n-j}{j} LS(j,k) = \sum_{j=0}^{n} \left[ \binom{2n-1-j}{j-1} + \binom{2n-1-j}{j} \right] LS(j,k)$$

$$= \sum_{j=0}^{n-1} \binom{2(n-1)-j}{j} LS(j+1,k) + \sum_{j=0}^{n-1} \binom{2(n-1)+1-j}{j} LS(j,k)$$

$$= \sum_{j=0}^{n-1} \binom{2(n-1)-j}{j} LS(j,k-1) + k(k+1) \sum_{j=0}^{n-1} \binom{2(n-1)-j}{j} LS(j,k) + \sum_{j=0}^{n-1} \binom{2(n-1)+1-j}{j} LS(j,k)$$

$$= T(n,k) + k(k+1)T(n,k+1) + (k+1)T(n,k+1) = T(n+1,k+1).$$

In the same way we get

$$\sum_{j=0}^{n}\binom{2n+1-j}{j}LS(j,k) = \sum_{j=0}^{n}\left(\binom{2n-j}{j-1}+\binom{2n-j}{j}\right)LS(j,k)$$

$$= \sum_{j=0}^{n-1}\binom{2(n-1)+1-j}{j}LS(j+1,k) + \sum_{j=0}^{n}\binom{2n-j}{j}LS(j,k)$$

$$= \sum_{j=0}^{n-1}\binom{2(n-1)+1-j}{j}LS(j,k-1) + k(k+1)\sum_{j=0}^{n-1}\binom{2(n-1)+1-j}{j}LS(j,k) + T(n+1,k+1)$$

$$= kT(n,k) + k(k+1)^2 T(n,k+1) + T(n+1,k+1) = (k+1)T(n+1,k+1).$$

**Corollary 3.3**

$$\left([j\leq i]\binom{2i-j}{j}\right)(LS(i,j)) = (T(i+1,j+1)) \tag{3.16}$$

*and*

$$\left([j\leq i]\binom{2i+1-j}{j}\right)(LS(i,j)) = (T(i+1,j+1))([j=i](j+1)). \tag{3.17}$$

**Corollary 3.4**

$$\left([j\leq i]\binom{2i+1-j}{j}\right)\left(\left([j\leq i]\binom{2i-j}{j}\right)\right)^{-1} = (T(i+1,j+1))([i=j](j+1))((T(i+1,j+1)))^{-1}. \tag{3.18}$$

*and*

$$\left(\left([j\leq i]\binom{2i-j}{j}\right)\right)^{-1}\left([j\leq i]\binom{2i+1-j}{j}\right) = (LS(i,j))([i=j](j+1))((LS(i,j)))^{-1}. \tag{3.19}$$

**Proof**

$$\left([j\leq i]\binom{2i+1-j}{j}\right)\left(\left([j\leq i]\binom{2i-j}{j}\right)\right)^{-1}$$

$$= (T(i+1,j+1))([i=j](j+1))((LS(i,j)))^{-1}(LS(i,j))((T(i+1,j+1)))^{-1}$$

In the same way we get (3.19).

Another interesting result is

**Theorem 3.5**

*The following identities hold:*

$$\sum_{j=0}^{n}\binom{n+1}{2n-2j}T(j+1,k+1) = LS(n+1,k+1) \quad (3.20)$$

*and*

$$\sum_{j=0}^{n}\binom{n+1}{2n-2j+1}T(j+1,k+1) = (k+1)LS(n+1,k+1). \quad (3.21)$$

*Equivalently this means*

$$\left([j \leq i]\binom{i+1}{2i-2j}\right)_{i,j\geq 0}\left(T(i+1,j+1)\right)_{i,j\geq 0} = \left(LS(i+1,j+1)\right)_{i,j\geq 0} \quad (3.22)$$

*and*

$$\left([j \leq i]\binom{i+1}{2i-2j+1}\right)_{i,j\geq 0}\left(T(i+1,j+1)\right)_{i,j\geq 0} = \left(LS(i+1,j+1)\right)_{i,j\geq 0}\left([i=j](j+1)\right)_{i,j\geq 0}. \quad (3.23)$$

**Proof**

For $n < k$ both sides vanish. For $n = k$ we get 1 and $k+1$ respectively.

$$\sum_{j=0}^{n}\binom{n+1}{2n-2j}T(j+1,k+1) = \sum_{j=0}^{n}\binom{n}{2n-2j}T(j+1,k+1) + \sum_{j=0}^{n}\binom{n}{2n-2j-1}T(j+1,k+1)$$

$$= \sum_{j=1}^{n}\binom{n}{2(n-1)+2-2j}T(j+1,k+1) + \sum_{j=0}^{n}\binom{n}{2(n-1)+1-2j}T(j+1,k+1)$$

$$= \sum_{j=0}^{n-1}\binom{n}{2(n-1)-2j}T(j+2,k+1) + (k+1)LS(n,k+1)$$

$$= \sum_{j=0}^{n-1}\binom{n}{2(n-1)-2j}T(j+1,k) + (k+1)^2\sum_{j=0}^{n-1}\binom{n}{2(n-1)-2j}T(j+1,k+1) + (k+1)LS(n,k+1)$$

$$= LS(n,k) + ((k+1)^2 + (k+1))LS(n,k+1) = LS(n+1,k+1).$$

For the second sum we get

$$\sum_{j=0}^{n}\binom{n+1}{2n-2j+1}T(j+1,k+1) = \sum_{j=0}^{n}\binom{n}{2n-2j+1}T(j+1,k+1) + \sum_{j=0}^{n}\binom{n}{2n-2j}T(j+1,k+1)$$

$$= \sum_{j=1}^{n}\binom{n}{2(n-1)+3-2j}T(j+1,k+1) + \sum_{j=1}^{n}\binom{n}{2(n-1)+2-2j}T(j+1,k+1)$$

$$= \sum_{j=0}^{n-1}\binom{n}{2(n-1)+1-2j}T(j+2,k+1) + \sum_{j=0}^{n-1}\binom{n}{2(n-1)-2j}T(j+2,k+1)$$

$$= \sum_{j=0}^{n-1}\binom{n}{2(n-1)+1-2j}T(j+1,k) + (k+1)^2\sum_{j=0}^{n-1}\binom{n}{2(n-1)+1-2j}T(j+1,k+1)$$

$$+ \sum_{j=0}^{n-1}\binom{n}{2(n-1)-2j}T(j+1,k) + (k+1)^2\sum_{j=0}^{n-1}\binom{n}{2(n-1)-2j}T(j+1,k+1)$$

$$= kLS(n,k) + (k+1)^2(k+1)LS(n,k+1) + LS(n,k) + (k+1)^2 LS(n,k+1)$$

$$= (k+1)\big(LS(n,k) + (k+1)(k+2)LS(n,k+1)\big) = (k+1)LS(n+1,k+1).$$

**Corollary 3.6**

$$\left(\left([j\le i]\binom{i+1}{2i-2j}\right)_{i,j\ge 0}\right)^{-1}\left([j\le i]\binom{i+1}{2i-2j+1}\right)_{i,j\ge 0} \quad (3.24)$$
$$= \big(T(i+1,j+1)\big)_{i,j\ge 0}\big([i=j](j+1)\big)_{i,j\ge 0}\left(\big(T(i+1,j+1)\big)_{i,j\ge 0}\right)^{-1}$$

and

$$\left([j\le i]\binom{i+1}{2i-2j+1}\right)_{i,j\ge 0}\left(\left([j\le i]\binom{i+1}{2i-2j}\right)_{i,j\ge 0}\right)^{-1} \quad (3.25)$$
$$= \big(LS(i+1,j+1)\big)_{i,j\ge 0}\big([i=j](j+1)\big)_{i,j\ge 0}\left(\big(LS(i+1,j+1)\big)_{i,j\ge 0}\right)^{-1}.$$

**Proof**

The left-hand side of (3.24) equals

$$\big(LS(i+1,j+1)\big)_{i,j\ge 0}\big([i=j](j+1)\big)_{i,j\ge 0}\left(\big(T(i+1,j+1)\big)_{i,j\ge 0}\right)^{-1}\big(T(i+1,j+1)\big)_{i,j\ge 0}\left(\big(LS(i+1,j+1)\big)_{i,j\ge 0}\right)^{-1}.$$

In the same way the left-hand side of (3.25) equals

$$\left(T(i+1,j+1)\right)_{i,j\geq 0}\left(\left(LS(i+1,j+1)\right)_{i,j\geq 0}\right)^{-1}\left(LS(i+1,j+1)\right)_{i,j\geq 0}\left([i=j](j+1)\right)_{i,j\geq 0}\left(\left(T(i+1,j+1)\right)_{i,j\geq 0}\right)^{-1}.$$

Comparing (3.24) with (3.18) we see that

$$\left([j\leq i]\binom{2i+1-j}{j}\right)\left(\left([j\leq i]\binom{2i-j}{j}\right)\right)^{-1}=\left(\left([j\leq i]\binom{i+1}{2i-2j}\right)_{i,j\geq 0}\right)^{-1}\left([j\leq i]\binom{i+1}{2i-2j+1}\right)_{i,j\geq 0}. \quad (3.26)$$

This implies

$$\left(\left([j\leq i]\binom{i+1}{2i-2j}\right)_{i,j\geq 0}\right)\left([j\leq i]\binom{2i+1-j}{j}\right)=\left([j\leq i]\binom{i+1}{2i-2j+1}\right)_{i,j\geq 0}\left(\left([j\leq i]\binom{2i-j}{j}\right)\right) \quad (3.27)$$
$$=\left(LS(i+1,j+1)\right)_{i,j\geq 0}\left([i=j](j+1)\right)_{i,j\geq 0}\left(\left(LS(i,j)\right)_{i,j\geq 0}\right)^{-1}.$$

**Remark**

Michael Schlosser has shown me a simple direct proof of the first identity in (3.27).

It suffices to show that

$$\sum_{j=k}^{n}\binom{n+1}{2n-2j}\binom{2j+1-k}{k}=\sum_{j=k}^{n}\binom{n+1}{2n-2j+1}\binom{2j-k}{k}.$$

This is equivalent with

$$\sum_{j=2k-1}^{2n}(-1)^{j+1}\binom{n+1}{2n-j}\binom{j+1-k}{k}=0.$$

This follows from the Chu-Vandermonde formula:

$$\sum_{j=2k-1}^{2n}(-1)^{j+1}\binom{n+1}{2n-j}\binom{j+1-k}{j+1-2k}$$
$$=\sum_{i=0}^{2n+1-2k}(-1)^{i}\binom{n+1}{2n+1-2k-i}\binom{i+k}{i}=\sum_{i=0}^{2+1-2k}\binom{n+1}{2n+1-2k-i}\binom{-k-1}{i}=\binom{n-k}{2n-2k+1}=0.$$

**3.5.** We later need the following result.

**Lemma 3.7**

*Let $w(0) = 1$ and $\hat{w}(n) = w(n+1)$. Then*

$$S^{\hat{w}}(n,k) = S^w(n+1, k+1) - S^w(n, k+1) \qquad (3.28)$$

*and*

$$s^{\hat{w}}(n,k) = -\sum_{j=0}^{k} s^w(n+1, j). \qquad (3.29)$$

*If*

$$F_1(n, j) = \sum_{\ell=0}^{n} S^w(n, \ell) F(\ell) s^w(\ell, j) = F(0) + \sum_{\ell=1}^{n} S^w(n, \ell) F(\ell) s^w(\ell, j) \qquad (3.30)$$

*for some function $F(\ell)$ then*

$$\sum_{\ell=0}^{n} S^{\hat{w}}(n, \ell) F(\ell+1) s^{\hat{w}}(\ell, k) = \sum_{j=0}^{k} \left( F_1(n, j) - F_1(n+1, j) \right). \qquad (3.31)$$

**Proof**

The first assertions follow from

$$\sum_n S^{\hat{w}}(n,k) x^n = \frac{x^k}{(1-w(1)x)(1-w(2)x)\cdots(1-w(k+1)x)} = \frac{(1-x)}{x} \frac{x^{k+1}}{(1-w(0)x)(1-w(1)x)\cdots(1-w(k+1)x)}$$

$$= (1-x) \sum_n S^w(n, k+1) x^{n-1} = \sum_n \left( S^w(n+1, k+1) - S^w(n, k+1) \right) x^n$$

and

$$\sum_{k=0}^{n} s^w(n,k) x^k = (x - w(0)) \cdots (x - w(n-1)) = (x-1) \sum_{k=0}^{n-1} s^{\hat{w}}(n-1, k) x^k$$

or

$$\sum_{k=0}^{n} s^{\hat{w}}(n,k) x^k = -\frac{1}{1-x} \sum_{k=0}^{n+1} s^w(n+1, k) x^k.$$

This implies

$$\sum_{\ell=0}^{n} S^{\hat{w}}(n,\ell)F(\ell+1)s^{\hat{w}}(\ell,k) = -\sum_{\ell=0}^{n}\left(S^{w}(n+1,\ell+1) - S^{w}(n,\ell+1)\right)F(\ell+1)\sum_{j=0}^{k}s^{w}(\ell+1,j)$$

$$= -\sum_{j=0}^{k}\sum_{\ell=0}^{n}S^{w}(n+1,\ell+1)F(\ell+1)s^{w}(\ell+1,j) + \sum_{j=0}^{k}\sum_{\ell=0}^{n-1}S^{w}(n,\ell+1)F(\ell+1)s^{w}(\ell+1,j)$$

$$= -\sum_{j=0}^{k}\sum_{\ell=1}^{n+1}S^{w}(n+1,\ell)F(\ell)s^{w}(\ell,j) + \sum_{j=0}^{k}\sum_{\ell=1}^{n}S^{w}(n,\ell)F(\ell)s^{w}(\ell,j)$$

$$= \sum_{j=0}^{k}\left(F_1(n,j) - F_1(n+1,j)\right).$$

## 4. Some interesting matrices

**4.1.** We know already that

$$\sum_{n\geq 0}\frac{F_n(s)}{n!}z^n = -e^z\sum_{n\geq 0}\frac{F_n(s)}{n!}(-z)^n. \tag{4.1}$$

This can also be written in the form

$$\sum\frac{F_{2n}(s)}{(2n)!}z^{2n} = \frac{e^z-1}{e^z+1}\sum\frac{F_{2n+1}(s)}{(2n+1)!}z^{2n+1}. \tag{4.2}$$

If we define a linear functional $\lambda$ on the vector space of polynomials in $s$ by

$$\lambda\left(F_{2n+1}(s)\right) = [n=0], \tag{4.3}$$

we get

$$\lambda(F_{2n}(s)) = \left[\frac{z^{2n}}{(2n)!}\right]z\frac{e^z-1}{e^z+1}. \tag{4.4}$$

Using (1.8) this gives (cf. [5],[3])

$$\lambda(F_{2n}(s)) = (-1)^{n-1}G_{2n}. \tag{4.5}$$

By comparing coefficients we get again (cf. [3])

$$F_{2n+2}(s) = \sum_{k=0}^{n}a(n,k)F_{2k+1}(s) \tag{4.6}$$

with

$$a(n,k) = (-1)^{n-k} \frac{1}{2k+1} \binom{2n+2}{2k} G_{2n-2k+2}. \tag{4.7}$$

We call the infinite matrix $A = (a(i,j))_{i,j \geq 0}$ and the finite parts $A_n = (a(i,j))_{i,j=0}^{n-1}$ **Genocchi matrices**.

E.g.

$$A_5 = \begin{pmatrix} 1 & 0 & 0 & 0 & 0 \\ -1 & 2 & 0 & 0 & 0 \\ 3 & -5 & 3 & 0 & 0 \\ -17 & 28 & -14 & 4 & 0 \\ 155 & -255 & 126 & -30 & 5 \end{pmatrix}. \tag{4.8}$$

It is clear that $A$ is the matrix version of "multiplication with $\dfrac{e^z - 1}{e^z + 1}$" of a certain kind of generating functions. More precisely we have

**Proposition 4.1**

Let $x = \begin{pmatrix} x_0 \\ x_1 \\ x_2 \\ \vdots \\ \vdots \end{pmatrix}$, $y = \begin{pmatrix} y_0 \\ y_2 \\ y_2 \\ \vdots \\ \vdots \end{pmatrix}$ and $A = (a(i,j))_{i,j \geq 0}$.

*Then $y = Ax$ is equivalent with*

$$\sum_{n \geq 0} y_n \frac{z^{2n+2}}{(2n+2)!} = \frac{e^z - 1}{e^z + 1} \sum_{n \geq 0} \frac{x_n}{(2n+1)!} z^{2n+1}. \tag{4.9}$$

For later uses we note that

$$\sum_{k=0}^{n} a(n,k) = 1. \tag{4.10}$$

This follows from (4.6) for $s = 0$.

For special values of $s$ (4.6) gives some interesting identities. Consider for example $s = -1$. Here $(F_n(s))_{n \geq 0} = (0,1,1,0,-1,-1,\cdots)$ is periodic with period 6.

This gives for example $\sum_{k=0}^{n} a(3n+2, 3k) = \sum_{k=0}^{n} a(3n+2, 3k+2)$.

For $s = -\frac{1}{4}$ we get $F_n\left(-\frac{1}{4}\right) = \frac{n}{2^{n-1}}$. Therefore

$$\sum_{k=0}^{n} 4^{n-k}(2k+1)a(n,k) = n+1.$$

By comparing coefficients of $s^k$ we get

$$\sum_{j=k}^{n} a(n, j) \binom{2j-k}{k} = \binom{2n+1-k}{k} \text{ for } n \geq k.$$

This is equivalent with the matrix identity

$$A = \left([j \leq i]\binom{2i+1-j}{j}\right)\left([j \leq i]\binom{2i-j}{j}\right)^{-1}. \tag{4.11}$$

If we recall that $\left([j \leq i]\binom{2i-j}{j}\right)_{i,j \geq 0} \begin{pmatrix} 1 \\ s \\ s^2 \\ s^3 \\ \vdots \end{pmatrix} = \begin{pmatrix} F_1 \\ F_3 \\ F_5 \\ F_7 \\ \vdots \end{pmatrix}$ and $\left([j \leq i]\binom{2i+1-j}{j}\right)_{i,j \geq 0} \begin{pmatrix} 1 \\ s \\ s^2 \\ s^3 \\ \vdots \end{pmatrix} = \begin{pmatrix} F_2 \\ F_4 \\ F_6 \\ F_8 \\ \vdots \end{pmatrix}$

we see that (4.11) is the same as (4.6).

By (3.18) we get the following representation

$$A = \left(T(i+1, j+1)\right)\left([i = j](j+1)\right)\left(\left(T(i+1, j+1)\right)\right)^{-1}. \tag{4.12}$$

By (3.26) we also have

$$A = \left(\left([j \leq i]\binom{i+1}{2i-2j}\right)_{i,j \geq 0}\right)^{-1} \left([j \leq i]\binom{i+1}{2i-2j+1}\right)_{i,j \geq 0}. \tag{4.13}$$

Thus we get

**Theorem 4.2**

*The Genocchi matrix A has the following representations*

$$A = \left([j \leq i]\binom{2i+1-j}{j}\right)\left([j \leq i]\binom{2i-j}{j}\right)^{-1}, \tag{4.14}$$

$$A = \left(\left([j \leq i]\binom{i+1}{2i-2j}\right)_{i,j \geq 0}\right)^{-1}\left([j \leq i]\binom{i+1}{2i-2j+1}\right)_{i,j \geq 0} \tag{4.15}$$

*and*

$$A = (T(i+1, j+1))_{i,j \geq 0}\,([i = j](i+1))_{i,j \geq 0}\,(t(i+1, j+1))_{i,j \geq 0}, \tag{4.16}$$

*where $t(i, j)$ are the central factorial numbers of the first kind and $T(i, j)$ the central factorial numbers of the second kind.*

Note the analogy with (3.9) and (3.10).

(4.15) implies $\left(\left([j \leq i]\binom{i+1}{2i-2j}\right)_{i,j \geq 0}\right) A = \left([j \leq i]\binom{i+1}{2i-2j+1}\right)_{i,j \geq 0}$.

By considering the first column of these matrices we get

$$\sum_{j=0}^{n}\binom{n+1}{2n-2j}(-1)^j G_{2j+2} = [n = 0].$$

Slightly reformulated this is

**Seidel's identity for Genocchi numbers ([13])**

$$\sum_{k=0}^{n}\binom{n}{2k}(-1)^k G_{2n-2k} = [n = 1]. \tag{4.17}$$

The matrix $\left([j \leq i]\binom{2i-j}{j}\right)^{-1}$ has been evaluated by Dumont and Zeng [6]. We don't need this result. We are only interested in the first column. Let $(H_{2n+1})_{n \geq 0} = (1, 1, 2, 8, 56, 608, \cdots)$ be the median Genocchi numbers (cf. OEIS A005439).

**Lemma 4.3**

*The elements of the first column of* $\left( [j \le i] \binom{2i-j}{j} \right)^{-1}$ *are the numbers* $(-1)^n H_{2n+1}$.

*The first two columns of*

$$\left( [j \le i] \binom{2i-j}{j} \right)^{-1} \left( [j \le i] \binom{2i+1-j}{j} \right)$$

*are* $\begin{pmatrix} 1 & 0 \\ 0 & 2 \\ 0 & -2 \\ 0 & 8 \\ \vdots & \vdots \end{pmatrix}$. *The numbers of the second column are*

$r(0) = -H_1 + 1 = 0, r(1) = H_3 + 1 = 2, r(2) = -H_5 = -2, r(3) = H_7 = 8, r(4) = -H_9 = -56, \cdots$ *and in general for* $n \ge 2$ $\quad r(n) = -\lambda(s^n) = (-1)^{n-1} H_{2n+1}$.

For example

$$\left( [j \le i] \binom{2i-j}{j} \right)^{-1}_{0 \le i,j \le 5} = \begin{pmatrix} 1 & 0 & 0 & 0 & 0 & 0 \\ -1 & 1 & 0 & 0 & 0 & 0 \\ 2 & -3 & 1 & 0 & 0 & 0 \\ -8 & 13 & -6 & 1 & 0 & 0 \\ 56 & -92 & 45 & -10 & 1 & 0 \\ -608 & 1000 & -493 & 115 & -15 & 1 \end{pmatrix}$$

$$\left( [j \le i] \binom{2i-j}{j} \right)^{-1}_{0 \le i,j \le 5} \left( [j \le i] \binom{2i+1-j}{j} \right)_{0 \le i,j \le 5} = \begin{pmatrix} 1 & 0 & 0 & 0 & 0 & 0 \\ 0 & 2 & 0 & 0 & 0 & 0 \\ 0 & -2 & 3 & 0 & 0 & 0 \\ 0 & 8 & -8 & 4 & 0 & 0 \\ 0 & -56 & 56 & -20 & 5 & 0 \\ 0 & 608 & -608 & 216 & -40 & 6 \end{pmatrix}.$$

To prove this lemma we need the fact (cf. [3],[5]) – which can serve as definition of the median Genocchi numbers – that

$$G_{2n} = \sum_{j=0}^{n} (-1)^{n-j} \binom{2n+1-j}{j} H_{2j+1}. \tag{4.18}$$

Since $\sum_{j=0}^{n+1}\binom{2n+1-j}{j}\lambda(s^j) = \lambda(F_{2n+2}(s)) = (-1)^n G_{2n+2}$ by (4.5) we see that

$$\lambda(s^n) = (-1)^n H_{2n+1}. \tag{4.19}$$

By (4.6) we get

$$\left([j \leq i]\binom{2i-j}{j}\right)^{-1}\begin{pmatrix}F_1(s)\\F_3(s)\\F_5(s)\\\vdots\end{pmatrix} = \begin{pmatrix}1\\s\\s^2\\\vdots\end{pmatrix} \tag{4.20}$$

and therefore by applying $\lambda$

$$\left([j \leq i]\binom{2i-j}{j}\right)^{-1}\begin{pmatrix}1\\0\\0\\\vdots\end{pmatrix} = \begin{pmatrix}1\\\lambda(s)\\\lambda(s^2)\\\vdots\end{pmatrix}.$$

To prove the second assertion we note that $F_{2n+1}(0) = 1$ and therefore (4.20) gives the first column.

To obtain the second column we recall that $\sum_{j=0}^{n}\binom{2n-j}{j}(-1)^j H_{2j+1} = \lambda(F_{2n+1}(s)) = 0$

and therefore $\sum_{j=0}^{n}\binom{2n-j}{j}r(j) = 1 + (2n-1) = 2n.$

(3.19) gives

$$\left(\left([j \leq i]\binom{2i-j}{j}\right)\right)^{-1}\left([j \leq i]\binom{2i+1-j}{j}\right) = (LS(i,j))([i=j](j+1))((LS(i,j)))^{-1}.$$

If we delete the first row and first column we get the matrix

$$(LS(i+1,j+1))_{i,j\geq 0}([i=j](j+2))_{i,j\geq 0}\left(((LS(i+1,j+1)))_{i,j\geq 0}\right)^{-1}. \tag{4.21}$$

**Theorem 4.4**

*Define linear functionals $\varphi_k$ by $\varphi_k(F_{2n-1}(s)) = T(n,k)$ for $k \geq 1$.*

*Then the following formulae hold:*

$$\varphi_{k+1}(s^n) = LS(n,k), \qquad (4.22)$$

$$\varphi_{k+1}(F_{2n+1}(s)) = \sum_{j=0}^{n-1} \binom{2n-j}{j} LS(j,k) = T(n+1,k+1) \qquad (4.23)$$

$$\varphi_{k+1}(F_{2n+2}(s)) = \sum_{j=0}^{n-1} \binom{2n+1-j}{j} LS(j,k) = (k+1)T(n+1,k+1). \qquad (4.24)$$

**Proof**

This is an immediate consequence of Theorem 3.2.

There is a simple Seidel-array for computing $LS(n,k)$ in terms of $T(n+1,k+1)$.

**Theorem 4.5**

*Let*

$$\begin{aligned}
&h(2i,0,k) = T(i+1,k+1), \\
&h(2i+1,0,k) = (k+1)T(i+1,k+1), \\
&h(i,j,k) = h(i,j-1,k) - h(i-1,j-1,k) \text{ if } j \leq \left\lfloor \frac{i}{2} \right\rfloor, \qquad (4.25) \\
&h(i,j,k) = 0 \text{ if } j > \left\lfloor \frac{i}{2} \right\rfloor.
\end{aligned}$$

*Then*

$$h(2n,n,k) = h(2n+1,n,k) = SL(n,k). \qquad (4.26)$$

For example for $k = 2$ this array begins with

$$\begin{pmatrix} \mathbf{0} & 0 & 0 & 0 & 0 & 0 & 0 & 0 & 0 \\ 0 & 0 & 0 & 0 & 0 & 0 & 0 & 0 & 0 \\ 0 & \mathbf{0} & 0 & 0 & 0 & 0 & 0 & 0 & 0 \\ 0 & 0 & 0 & 0 & 0 & 0 & 0 & 0 & 0 \\ 1 & 1 & \mathbf{1} & 0 & 0 & 0 & 0 & 0 & 0 \\ 3 & 2 & 1 & 0 & 0 & 0 & 0 & 0 & 0 \\ 14 & 11 & 9 & \mathbf{8} & 0 & 0 & 0 & 0 & 0 \\ 42 & 28 & 17 & 8 & 0 & 0 & 0 & 0 & 0 \\ 147 & 105 & 77 & 60 & \mathbf{52} & 0 & 0 & 0 & 0 \end{pmatrix}.$$

**Proof**

It suffices to show that $h(2i, j, k) = \varphi_k\left(F_{2i-2j+1}(s)s^j\right)$ and $h(2i+1, j, k) = \varphi_k\left(F_{2i-2j+2}(s)s^j\right)$.

For then we have $h(2i, i, k) = \varphi_{k+1}\left(F_1(s)s^i\right) = \varphi_{k+1}(s^i) = LS(i,k)$ and

$h(2i+1, i, k) = \varphi_{k+1}\left(F_2(s)s^i\right) = \varphi_{k+1}(s^i) = LS(i,k)$.

For $j = 0$ this is clear by our definition.

We have only to verify that $h(2i, j, k) = h(2i, j-1, k) - h(2i-1, j-1, k)$

and $h(2i+1, j, k) = h(2i+1, j-1, k) - h(2i, j-1, k)$ for $j \leq i$.

But this follows from

$F_{2i-2j+1}(s)s^j = F_{2i-2j+3}(s)s^{j-1} - F_{2i-2j+2}(s)s^{j-1} = F_{2i-2(j-1)+1}(s)s^{j-1} - F_{2i-2(j-1)}(s)s^{j-1}$

and

$F_{2i-2j+2}(s)s^j = F_{2i-2j+4}(s)s^{j-1} - F_{2i-2j+3}(s)s^{j-1} = F_{2i-2(j-1)+2}(s)s^{j-1} - F_{2i-2(j-1)+1}(s)s^{j-1}$.

Note that $h(2n+1, 0, k) = \sum_{j=0}^{n} h(2n, j, k)$. This follows from the identity

$F_{2n+2}(s) = \sum_{j=0}^{n} F_{2n+1-2j}(s)s^j$.

**Remark**

If we replace $\varphi_k\left(F_{2i-2j+1}(s)s^j\right)$ by $h(2i, j) = \lambda\left(F_{2i-2j+1}(s)s^j\right)$ we get the original Seidel triangle for the Genocchi numbers in [13], Beilage 2, which Seidel called "Differenzen-Tableau der Bernoulli'schen Zähler".

Its first terms are

$$\begin{pmatrix}
1 & 0 & 0 & 0 & 0 & 0 & 0 & 0 & 0 & 0 \\
1 & 0 & 0 & 0 & 0 & 0 & 0 & 0 & 0 & 0 \\
0 & -1 & 0 & 0 & 0 & 0 & 0 & 0 & 0 & 0 \\
-1 & -1 & 0 & 0 & 0 & 0 & 0 & 0 & 0 & 0 \\
0 & 1 & 2 & 0 & 0 & 0 & 0 & 0 & 0 & 0 \\
3 & 3 & 2 & 0 & 0 & 0 & 0 & 0 & 0 & 0 \\
0 & -3 & -6 & -8 & 0 & 0 & 0 & 0 & 0 & 0 \\
-17 & -17 & -14 & -8 & 0 & 0 & 0 & 0 & 0 & 0 \\
0 & 17 & 34 & 48 & 56 & 0 & 0 & 0 & 0 & 0 \\
155 & 155 & 138 & 104 & 56 & 0 & 0 & 0 & 0 & 0
\end{pmatrix}.$$

Here we have $h(2n,n) = (-1)^n H_{2n+1}$.

It has the special property that we can also compute the first column without previous knowledge of the Genocchi numbers because $h(2n+1, 0) = \sum_{j=0}^{n} h(2n, j)$ and the fact that

$h(2n, 0) = \lambda(F_{2n+1}(s)) = 0$ for $n > 0$.

**4.2.** For our next results we need the square of $A$.

For example

$$A_5^2 = \begin{pmatrix}
1 & 0 & 0 & 0 & 0 \\
-3 & 4 & 0 & 0 & 0 \\
17 & -25 & 9 & 0 & 0 \\
-155 & 238 & -98 & 16 & 0 \\
2073 & -3255 & 1428 & -270 & 25
\end{pmatrix}.$$

**Theorem 4.6**

*The square of $A_n$ is given by*

$$A_n^2 = \left(aa(i,j)[j \leq i]\right)_{i,j=0}^{n-1} \tag{4.27}$$

*with*

$$aa(n,k) = (-1)^{n-k} \binom{2n+2}{2k} \frac{n+k+2}{(2k+1)(n+2-k)} G_{2n-2k+4}. \tag{4.28}$$

**Proof**

Let $f(z) = \sum_{n \geq 0} \frac{x_n}{(2n+1)!} z^{2n+1}$.

Since $y = Ax$ is equivalent with $\sum_{n \geq 0} y_n \frac{z^{2n+2}}{(2n+2)!} = \frac{e^z - 1}{e^z + 1} f(z)$ we see that $z = Ay = A^2 x$ is equivalent with

$$\sum_{n \geq 1} z_{n-1} \frac{z^{2n}}{(2n)!} = \frac{e^z - 1}{e^z + 1} \sum_{n \geq 1} y_{n-1} \frac{z^{2n-1}}{(2n-1)!} = \frac{e^z - 1}{e^z + 1} \left( \frac{e^z - 1}{e^z + 1} f(z) \right)'$$

$$= \left( \frac{e^z - 1}{e^z + 1} \right)^2 f'(z) + \frac{e^z - 1}{e^z + 1} \left( \frac{e^z - 1}{e^z + 1} \right)' f(z) \qquad (4.29)$$

Since

$$\left( \frac{e^z - 1}{e^z + 1} \right)^2 = \sum_{n \geq 1} (-1)^{n+1} \frac{G_{2n+2}}{n+1} \frac{z^{2n}}{(2n)!} \qquad (4.30)$$

and

$$-\frac{e^z - 1}{e^z + 1} \left( \frac{e^z - 1}{e^z + 1} \right)' = \left( \frac{e^z - 1}{e^z + 1} \right)'' = \sum_{n \geq 1} (-1)^{n-1} \frac{G_{2n}}{2n} \frac{z^{2n-3}}{(2n-3)!}, \qquad (4.31)$$

we get

$$\left( \frac{e^z - 1}{e^z + 1} \right)^2 f'(z) = \sum_n \frac{z^{2n}}{(2n)!} \sum_{k=1}^{n} \binom{2n}{2k-2} (-1)^{n-k} \frac{G_{2n-2k+4}}{n-k+2} x(k) \qquad (4.32)$$

and

$$-\left( \frac{e^z - 1}{e^z + 1} \right)'' f(z) = \sum_n \frac{z^{2n}}{(2n)!} \sum_{k=0}^{n} \binom{2n}{2k-1} (-1)^{n-k} \frac{G_{2n-2k+4}}{2n-2k+4} x(k). \qquad (4.33)$$

From (4.32) and (4.33) and using (4.29) we get immediately (4.28).

It only remains to prove (4.30).

This follows from

$$\left(\frac{e^z-1}{e^z+1}\right)^2 = 1-2\left(\frac{e^z-1}{e^z+1}\right)' = 1-2\sum_{n\geq 0}(-1)^n G_{2n+2}\frac{(2n+1)z^{2n}}{(2n)!(2n+1)(2n+2)}$$

$$=\sum_{n\geq 0}(-1)^{n-1}\frac{G_{2n+2}}{n+1}\frac{z^{2n}}{(2n)!}.$$

**Theorem 4.7 (Further properties of $A_n^2$)**

$$aa(n,0) = -a(n+1,0) \tag{4.34}$$

*and for $1 \leq k \leq n-1$*

$$aa(n,k) = a(n,k-1) - a(n+1,k) \tag{4.35}$$

**Proof.**

This follows from (4.28).

**4.3.** Intimately connected with the linear functional $\lambda$ is the linear functional $\lambda^*$ defined by

$$\lambda^*(F_n(s)) = -\lambda(sF_n(s)). \tag{4.36}$$

Since $\lambda(sF_n(s)) = \lambda(F_{n+2}(s)) - \lambda(F_{n+1}(s))$ we get

$\lambda(sF_{2n}(s)) = \lambda(F_{2n+2}(s)) - \lambda(F_{2n+1}(s)) = \lambda(F_{2n+2}(s))$ for $n>0$ and

$\lambda(sF_{2n+1}(s)) = \lambda(F_{2n+3}(s)) - \lambda(F_{2n+2}(s)) = -\lambda(F_{2n+2}(s)).$

The sequence $\left(\lambda^*(F_n(s))\right)$ begins with

$0,1,1,-1,-3,3,17,-17,-155,155,2073,-2073,-38227.$

Therefore

$$\lambda^*(F_{2n}(s) + F_{2n+1}(s)) = [n=0] \tag{4.37}$$

and

$$\lambda^*(F_{2n-1}(s) + F_{2n}(s)) = (-1)^n (G_{2n} + G_{2n+2}). \tag{4.38}$$

Thus we are led to consider the basis consisting of the polynomials $F_{2k-2}(s)+F_{2k-1}(s)$.

Here we get

**Theorem 4.8**

*Let*

$$a_1(n,k) = [k \le n]\sum_{j=0}^{k} a(n,j). \tag{4.39}$$

*Then*

$$F_{2n+1}(s) = \sum_{k=0}^{n} a_1(n,k)\bigl(F_{2k}(s)+F_{2k+1}(s)\bigr). \tag{4.40}$$

The first entries of $(a_1(n,k))$ are

$$\begin{pmatrix} 1 & 0 & 0 & 0 & 0 & 0 \\ -1 & 1 & 0 & 0 & 0 & 0 \\ 3 & -2 & 1 & 0 & 0 & 0 \\ -17 & 11 & -3 & 1 & 0 & 0 \\ 155 & -100 & 26 & -4 & 1 & 0 \\ -2073 & 1337 & -346 & 50 & -5 & 1 \end{pmatrix}.$$

**Proof**

We know that $F_{2n}(s) = \sum_{k=0}^{n-1} a(n-1,k)F_{2k+1}(s)$. Therefore

$F_{2n}(s)+F_{2n+1}(s) = \sum_{k=0}^{n-1} a(n-1,k)F_{2k+1}(s)+F_{2n+1}(s)$. This implies that in (4.40) the matrix $(a_1(i,j))$ is the inverse of $I+B$, where $I$ denotes the identity matrix and $b(n,k) = a(n-1,k)$ for $k<n$. We have to show that $a_1(n,k)$ is given by (4.39).

It suffices to show that $(I+B)\left([j \le i]\sum_{\ell=0}^{j} a(i,\ell)\right) = I.$

This means that for each $k$

$$\sum_{j=k}^{n-1} a(n-1,j)\bigl(a(j,0)+a(j,1)+\cdots a(j,k)\bigr)+a(n,0)+\cdots+a(n,k) = [n=k].$$

Let $k < n$. By definition $aa(n-1,k) = \sum_{j=0}^{n-1} a(n-1,j)a(j,k)$. Therefore

$$aa(n-1,i) - \sum_{j=k}^{n-1} a(n-1,j)a(j,i) = \sum_{j=0}^{k-1} a(n-1,j)a(j,i).$$

Using (4.35) and (4.10) the left-hand side is equivalent with

$$aa(n-1,0) - \sum_{j=0}^{k-1} a(n-1,j)a(j,0) + aa(n-1,1) - \sum_{j=1}^{k-1} a(n-1,j)a(j,1) + \cdots$$
$$+ aa(n-1,k-1) - a(n-1,k-1)a(k-1,k-1) + aa(n-1,k) + a(n,0) + a(n,1) + \cdots + a(n,k)$$
$$= aa(n-1,0) + a(n,0) + aa(n-1,1) + a(n,1) + \cdots + aa(n-1,k) + a(n,k) - \sum_{i=0}^{k-1} a(n-1,i) \sum_{j=0}^{i} a(i,j)$$
$$= a(n-1,0) + a(n-1,1) + \cdots + a(n-1,k-1) - \sum_{j=0}^{k-1} a(n-1,j) = 0.$$

For $k = n$ the first sum vanishes and $a(n,0) + \cdots + a(n,n) = 1$.

**Corollary 4.9**

*Let*

$$a_2(n,k) = [k \leq n](a_1(n,k) - a_1(n+1,k)). \tag{4.41}$$

*Then*

$$F_{2n+1}(s) + F_{2n+2}(s) = \sum_{k=0}^{n} a_2(n,k)(F_{2k}(s) + F_{2k+1}(s)). \tag{4.42}$$

**Proof**

By Theorem 4.7 we have

$$F_{2n+1}(s) = \sum_{k=0}^{n} a_1(n,k)(F_{2k}(s) + F_{2k+1}(s)) = \sum_{k=0}^{n-1} a_1(n,k)(F_{2k}(s) + F_{2k+1}(s)) + F_{2n}(s) + F_{2n+1}(s).$$

Therefore

$$F_{2n+2}(s) = -\sum_{k=0}^{n} a_1(n+1,k)(F_{2k}(s) + F_{2k+1}(s))$$

and thus $F_{2n+1}(s) + F_{2n+2}(s) = \sum_{k=0}^{n} a_2(n,k)(F_{2k}(s) + F_{2k+1}(s)).$

The first terms of $(a_2(n,k))$ are

$$\begin{pmatrix} 2 & 0 & 0 & 0 & 0 \\ -4 & 3 & 0 & 0 & 0 \\ 20 & -13 & 4 & 0 & 0 \\ -172 & 111 & -29 & 5 & 0 \\ 2228 & -1437 & 372 & -54 & 6 \end{pmatrix}.$$

Since $a(n,k) = \sum_{j=0}^{n} T(n+1, j+1)(j+1)t(j+1,k+1)$ we see that $a_2(n,k) = F(n, j)$ if in Lemma 3.11 we choose $w(n) = (n+1)^2$ and $F(\ell) = \ell + 1$.

Therefore $a_2(n,k) = \sum_{\ell=0}^{n} S^{\hat{w}}(n,\ell) F(\ell+1) s^{\hat{w}}(\ell,k)$.

This gives

**Theorem 4.10**

Let $w(n) = (n+1)^2$ and thus $\hat{w}(n) = (n+2)^2$. Then

$$\left(a_2(i,j)\right)_{i,j=0}^{n-1} = \left(S^{\hat{w}}(i,j)\right)_{i,j=0}^{n-1} \left([i=j](i+2)\right)_{i,j=0}^{n-1} \left(s^{\hat{w}}(i,j)\right)_{i,j=0}^{n-1}. \quad (4.43)$$

**4.5.** Next we consider the linear functional $\mu$ defined by $\mu(F_{2k+2}(s)) = [k = 0]$. From (4.2) we see that

$$\sum \frac{\mu(F_{2n+1}(s))}{(2n+1)!} z^{2n+1} = \frac{e^z + 1}{e^z - 1} \sum \frac{\mu(F_{2n}(s))}{(2n)!} z^{2n} = \frac{z^2}{2} \frac{e^z + 1}{e^z - 1} = \sum_n (2n+1) B_{2n} \frac{z^{2n+1}}{(2n+1)!}. \quad (4.44)$$

This implies

$$\mu\left(F_{2n+1}(s)\right) = (2n+1) B_{2n} \quad (4.45)$$

and

$$F_{2n+1}(s) = \sum_{j=0}^{n} \binom{2n+1}{2j+1} \frac{B_{2n-2j}}{j+1} F_{2j+2}(s). \quad (4.46)$$

Comparing (4.6) with (4.46) we see that

$$A^{-1} = (a(j,k))_{i,j \geq 0}^{-1} = \left( \binom{2j+1}{2k+1} \frac{B_{2j-2k}}{k+1} \right)_{i,j \geq 0}. \tag{4.47}$$

The first terms of the sequence $((2n+1)B_{2n})_{n \geq 0}$ are $1, \frac{1}{2}, -\frac{1}{6}, \frac{1}{6}, -\frac{3}{10}, \frac{5}{6}, -\frac{691}{210}, \frac{35}{2}, -\frac{3617}{30}$.

Let $w(j,k) = \binom{2j+1}{2k+1} \frac{B_{2j-2k}}{k+1} [k \leq j]$.

By (4.15) we have

$$\left( [j \leq i] \binom{i+1}{2i-2j+1} \right)_{i,j \geq 0} A^{-1} = \left( [j \leq i] \binom{i+1}{2i-2j} \right)_{i,j \geq 0}.$$

Considering the first column we get

$$\sum_{j=0}^{n} \binom{n+1}{2n-2j+1} (2j+1) B_{2j} = \sum_{j=0}^{n} \binom{n+1}{2j+1} (2n-2j+1) B_{2n-2j}$$
$$= \sum_{j=0}^{n} \binom{n+1}{n-2j} (n+n-2j+1) B_{n+n-2j} = [n=0].$$

Since $B_{2i+1} = 0$ for $i > 0$ this is equivalent with

$$\sum_{i=0}^{n+1} \binom{n+1}{i} (n+i+1) B_{n+i} = 0 \text{ for } n > 1. \text{ But it also holds for } n=0 \text{ and } n=1.$$

Therefore we get

**Kaneko's identity ([7],[10],[13])**

$$\sum_{i=0}^{n+1} \binom{n+1}{i} (n+i+1) B_{n+i} = 0. \tag{4.48}$$

This identity has first been proved by A.v. Ettingshausen [7] and has been rediscovered by L. Seidel [13],VIII, and by M. Kaneko [10].

(4.47) implies that

$$\left(w(j,k)\right)_{j,k=0}^{n-1} = \left(T(i+1, j+1)\right)_{i,j=0}^{n-1} \left([i=j]\frac{1}{j+1}\right)_{i,j=0}^{n-1} \left(t(i+1, j+1)\right)_{i,j=0}^{n-1}. \tag{4.49}$$

The inverse of $\left(a_2(i,j)\right)_{i,j=0}^{n-1} = \left(S^{\hat{w}}(i,j)\right)_{i,j=0}^{n-1} \left([i=j](i+2)\right)_{i,j=0}^{n-1} \left(s^{\hat{w}}(i,j)\right)_{i,j=0}^{n-1}$ is

$$\left(z(i,j)\right)_{i,j=0}^{n-1} = \left(S^{\hat{w}}(i,j)\right)_{i,j=0}^{n-1} \left([i=j]\frac{1}{i+2}\right)_{i,j=0}^{n-1} \left(s^{\hat{w}}(i,j)\right)_{i,j=0}^{n-1}.$$

By Lemma 3.11 this implies that $z(n,k) = \sum_{j=0}^{k} \left(w(n,j) - w(n+1,j)\right)$ for $k \leq n$ and $z(n,k) = 0$ else.

So we have

$$F_{2n}(s) + F_{2n+1}(s) = \sum_{k=0}^{n} z(n,k)\left(F_{2k+1}(s) + F_{2k+2}(s)\right). \tag{4.50}$$

## 5. Analogous results for Lucas polynomials

**5.1.** For the Lucas polynomials we get

$$\sum \frac{L_{2n+1}(s)}{(2n+1)!} z^{2n+1} = \frac{e^z - 1}{e^z + 1} \sum \frac{L_{2n}(s)}{(2n)!} z^{2n}. \tag{5.1}$$

This follows from

$$e^{\alpha z} + e^{\beta z} = e^{(1-\beta)z} + e^{(1-\alpha)z} = e^z \left(e^{-\alpha z} + e^{-\beta z}\right). \tag{5.2}$$

We write the series expansion of $\dfrac{e^z - 1}{e^z + 1}$ in the form

$$\frac{e^z - 1}{e^z + 1} = \sum_{k \geq 0} (-1)^k \frac{T_{2k+1}}{2^{2k+1}} \frac{z^{2k+1}}{(2k+1)!}, \tag{5.3}$$

where $T_n$, $n \geq 1$, are the tangent numbers $1, 2, 16, 272, 7936, \cdots$.

Therefore

$$\sum \frac{Y_n}{(2n+1)!} z^{2n+1} = \frac{e^z - 1}{e^z + 1} \sum \frac{X_n}{(2n)!} z^{2n} \quad (5.4)$$

is equivalent with

$$Y_n = \sum_{k=0}^{n} (-1)^k \frac{T_{2k+1}}{2^{2k+1}} \binom{2n+1}{2k+1} X_{n-k} = \sum_{k=0}^{n} (-1)^{n-k} \frac{T_{2n-2k+1}}{2^{2n-2k+1}} \binom{2n+1}{2k} X_k.$$

If we set $b(n,k) = (-1)^{n-k} \frac{T_{2n-2k+1}}{2^{2n-2k+1}} \binom{2n+1}{2k} = (-1)^{n-k} \frac{1}{2} \binom{2n+1}{2k} \frac{G_{2n-2k+2}}{n-k+1}$, then (5.4) is the same as

$$Y_n = \sum_{j=0}^{n} b(n, j) X_j. \quad (5.5)$$

We call the matrices

$$B = (b(i, j))_{i, j \geq 0} = \left( (-1)^{i-j} \frac{T_{2i-2j+1}}{2^{2i-2j+1}} \binom{2i+1}{2j} \right)_{i, j \geq 0} \quad (5.6)$$

**tangent-matrices**.

For example

$$B_5 = \begin{pmatrix} \frac{1}{2} & 0 & 0 & 0 & 0 \\ -\frac{1}{4} & \frac{3}{2} & 0 & 0 & 0 \\ \frac{1}{2} & -\frac{5}{2} & \frac{5}{2} & 0 & 0 \\ -\frac{17}{8} & \frac{21}{2} & -\frac{35}{4} & \frac{7}{2} & 0 \\ \frac{31}{2} & -\frac{153}{2} & 63 & -21 & \frac{9}{2} \end{pmatrix}.$$

If we compare (5.1) with (5.5) we see that $L_{2n+1}(s) = \sum_{j=0}^{n} b(n, j) L_{2j}(s)$. This is again (2.3).

Let $L_n(s) = \sum_{k=0}^{\lfloor \frac{n}{2} \rfloor} l(n, s) s^k$. Then $L_{2n}(s) = \sum_{k=0}^{n} l(2n, k) s^k$ and $L_{2n+1}(s) = \sum_{k=0}^{n} l(2n+1, k) s^k$.

Therefore we get

$$B = \left(l(2i+1, j)\right)_{i,j \geq 0} \left(\left(l(2i, j)\right)_{i,j \geq 0}\right)^{-1}. \tag{5.7}$$

**5.2.** Let $w(n) = \left(\dfrac{2n+1}{2}\right)^2$ and let $U(n,k) = S^w(n,k)$ and $u(n,k) = s^w(n,k)$. respectively.

The first values of the numbers $4^{n-k}U(n,k)$ and $4^{n-k}u(n,k)$ are given by the following tables.

$$\left(4^{i-j}U(i,j)\right)_{i,j=0}^{6} = \begin{pmatrix} 1 & 0 & 0 & 0 & 0 & 0 & 0 \\ 1 & 1 & 0 & 0 & 0 & 0 & 0 \\ 1 & 10 & 1 & 0 & 0 & 0 & 0 \\ 1 & 91 & 35 & 1 & 0 & 0 & 0 \\ 1 & 820 & 966 & 84 & 1 & 0 & 0 \\ 1 & 7381 & 24970 & 5082 & 165 & 1 & 0 \\ 1 & 66430 & 631631 & 273988 & 18447 & 286 & 1 \end{pmatrix}$$

and

$$\left(4^{i-j}u(i,j)\right)_{i,j=0}^{6} = \begin{pmatrix} 1 & 0 & 0 & 0 & 0 & 0 & 0 \\ -1 & 1 & 0 & 0 & 0 & 0 & 0 \\ 9 & -10 & 1 & 0 & 0 & 0 & 0 \\ -225 & 259 & -35 & 1 & 0 & 0 & 0 \\ 11025 & -12916 & 1974 & -84 & 1 & 0 & 0 \\ -893025 & 1057221 & -172810 & 8778 & -165 & 1 & 0 \\ 108056025 & -128816766 & 21967231 & -1234948 & 28743 & -286 & 1 \end{pmatrix}.$$

There is also an analogue of Theorem 3.2.

To this end we introduce an analogue of the Legendre-Stirling numbers.

Let $w(n) = \dfrac{(2n-1)(2n+1)}{4}$ and $V(n,k) = S^w(n,k)$ and $v(n,k) = s^w(n,k)$.

Then we get

**Theorem 5.1**

$$\sum_{j=0}^{n} l(2n, j)V(j,k) = 2U(n,k) \tag{5.8}$$

*and*

$$\sum_{j=0}^{n} l(2n+1, j)V(j,k) = (2k+1)U(n,k). \tag{5.9}$$

**Proof**

We use again induction. Let both identities be already proved for $n-1$. Then

$$\sum_{j=0}^{n} l(2n, j)V(j,k) = \sum_{j=0}^{n-1} l(2n-1, j)V(j,k) + \sum_{j=0}^{n-1} l(2n-2, j-1)V(j,k)$$

$$= \sum_{j=0}^{n-1} l(2(n-1)+1, j)V(j,k) + \sum_{j=0}^{n-1} l(2(n-1), j)V(j+1,k)$$

$$= (2k+1)U(n-1,k) + \sum_{j=0}^{n-1} l(2(n-1), j)V(j,k-1) + \frac{(2k-1)(2k+1)}{4}\sum_{j=0}^{n-1} l(2(n-1), j)V(j,k)$$

$$= (2k+1)U(n-1,k) + 2U(n-1,k-1) + \frac{(2k-1)(2k+1)}{4} 2U(n-1,k)$$

$$= 2\left(U(n-1,k-1) + \frac{(2k+1)^2}{4}U(n-1,k)\right) = 2U(n,k).$$

And

$$\sum_{j=0}^{n} l(2n+1, j)V(j,k) = \sum_{j=0}^{n} l(2n, j)V(j,k) + \sum_{j=0}^{n-1} l(2n-1, j-1)V(j,k)$$

$$= 2U(n,k) + \sum_{j=0}^{n-1} l(2n-1, j)V(j+1,k)$$

$$= 2U(n,k) + \sum_{j=0}^{n-1} l(2n-1, j)V(j,k-1) + \frac{(2k-1)(2k+1)}{4}\sum_{j=0}^{n-1} l(2(n-1), j)V(j,k)$$

$$= 2U(n,k) + (2k-1)U(n-1,k-1) + \frac{(2k-1)(2k+1)}{4}(2k+1)U(n-1,k)$$

$$= 2\left(U(n,k) + \frac{2k-1}{2}\left(U(n-1,k-1) + \frac{(2k+1)^2}{4}U(n-1,k)\right)\right)$$

$$= 2\left(U(n,k) + \frac{2k-1}{2}U(n,k)\right) = (2k+1)U(n,k).$$

This implies

**Theorem 5.2**

*The tangent matrix B has the factorization*

$$B = \left(l(2i+1, j)\right)_{i,j\geq 0}\left(\left(l(2i, j)\right)_{i,j\geq 0}\right)^{-1} = \left(U(i,j)\right)_{i,j\geq 0}\left([i=j]\frac{2j+1}{2}\right)_{i,j\geq 0}\left(u(i,j)\right)_{i,j\geq 0}, \quad (5.10)$$

where $u(i, j)$ and $U(i, j)$ are the generalized Stirling numbers corresponding to
$$w(n) = \left(\frac{2n+1}{2}\right)^2.$$

In the same way as in Theorem 4.5 there is a simple Seidel-array for computing $V(n,k)$ in terms of $U(n,k)$.

**Proposition 5.3**

*Let*

$$h(2i,0,k) = U(i,k),$$
$$h(2i+1,0,k) = \frac{2k+1}{2}U(i,k),$$
$$h(i,j,k) = h(i,j-1,k) - h(i-1,j-1,k) \text{ if } j \leq \left\lfloor \frac{i}{2} \right\rfloor, \tag{5.11}$$
$$h(i,j,k) = 0 \text{ if } j > \left\lfloor \frac{i}{2} \right\rfloor.$$

*Then*

$$h(2n,n,k) = V(n,k). \tag{5.12}$$

For example $(h(i,j,1))$ begins with

$$\begin{pmatrix} \mathbf{0} & 0 & 0 & 0 & 0 & 0 & 0 \\ 0 & 0 & 0 & 0 & 0 & 0 & 0 \\ 1 & \mathbf{1} & 0 & 0 & 0 & 0 & 0 \\ \frac{3}{2} & \frac{1}{2} & 0 & 0 & 0 & 0 & 0 \\ \frac{5}{2} & 1 & \frac{1}{2} & 0 & 0 & 0 & 0 \\ \frac{15}{4} & \frac{5}{4} & \frac{1}{4} & 0 & 0 & 0 & 0 \\ \frac{91}{16} & \frac{31}{16} & \frac{11}{16} & \frac{7}{16} & 0 & 0 & 0 \end{pmatrix}.$$

**Proposition 5.4**

*The inverse of $B$ is*

$$B^{-1} = \left([j \leq i] \binom{2i}{2j} \frac{B_{2i-2j}}{2j+1}\right). \tag{5.13}$$

This follows from (5.4) and the well-known series $\dfrac{z}{2}\dfrac{e^z+1}{e^z-1} = \sum_{n\geq 0} B_{2n}\dfrac{z^{2n}}{(2n)!}$.

Then from (5.1) we get again (2.4).

## 6. Some interesting identities

**6.1.** For any sequence $(a(n))_{n\geq 0}$ the sequence $c(n) = \sum_{j=0}^{n}(-1)^j\left(\prod_{i=0}^{j-1}w(i)\right)S^w(n,j)a(j)$ is the first column in $\left(S^w(i,j)\right)_{i,j=0}^{\infty}\left(a(i)([i=j])\right)_{i,j=0}^{\infty}\left(s^w(i,j)\right)_{i,j=0}^{\infty}$.

More precisely we have

**Theorem 6.1**

*Let*

$$X = \left(x(i,j)\right)_{i,j\geq 0} = \left(S^w(i,j)\right)_{i,j\geq 0}\left(a(i)([i=j])\right)_{i,j\geq 0}\left(s^w(i,j)\right)_{i,j\geq 0}. \quad (6.1)$$

*Then*

$$x(n,0) = \sum_{j=0}^{n}(-1)^j S^w(n,j)a(j)\left(\prod_{i=0}^{j-1}w(i)\right) \quad (6.2)$$

*and*

$$\sum_{j=0}^{n} s^w(n,j)x(j,0) = (-1)^n a(n)\prod_{j=0}^{n-1}w(j). \quad (6.3)$$

The sequence $(c(n))_{n\geq 0}$ can be simply computed with the Akiyama-Tanigawa algorithm (cf. [1], [9], [11], [15] ):

**Theorem 6.2 ( Akiyama-Tanigawa algorithm)  ( [1],[9],[11],[15])**

*Suppose that* $w(n) \neq 0$ *for all* $n$. *Let* $c(n) = \sum_{j=0}^{n}(-1)^j \left(\prod_{i=0}^{j-1} w(i)\right) S^w(n, j) a(j)$.

*Define a matrix* $M = (m(i,j))_{i,j \geq 0}$ *by*

$m(0, j) = a(j)$ *for* $j \in \mathbb{N}$ *and*

$$m(i, j) = w(j)\bigl(m(i-1, j) - m(i-1, j+1)\bigr) \qquad (6.4)$$

*Then* $m(n,0) = c(n) = \sum_{j=0}^{n}(-1)^j \left(\prod_{i=0}^{j-1} w(i)\right) S^w(n, j) a(j)$.

**Proof**

To prove this we show more generally that

$$m(n,k) = \frac{(-1)^k}{\prod_{\ell=0}^{k-1} w(\ell)} \sum_{i=0}^{k} s(k,i) c(n+i). \qquad (6.5)$$

This holds for $n = 0$ because

$$\frac{(-1)^k}{\prod_{\ell=0}^{k-1} w(\ell)} \sum_{i \geq 0} s^w(k,i) c(i) = \frac{(-1)^k}{\prod_{\ell=0}^{k-1} w(\ell)} \sum_{i \geq 0} s^w(k,i) \sum_{j=0}^{i}(-1)^j \left(\prod_{\ell=0}^{j-1} w(\ell)\right) S^w(i, j) a(j)$$

$$= \frac{(-1)^k}{\prod_{\ell=0}^{k-1} w(\ell)} \sum_{j \geq 0}(-1)^j \left(\prod_{\ell=0}^{j-1} w(\ell)\right) a(j) \sum_{i \geq 0} s^w(k,i) S^w(i, j) = a(k).$$

Now suppose that (6.5) holds for $n-1$. Then

$$w(k)\big(m(n-1,k)-m(n-1,k+1)\big)$$

$$= w(k)\left(\frac{(-1)^k}{\prod_{\ell=0}^{k-1}w(\ell)}\sum_{i\geq 0}s^w(k,i)c(n-1+i) - \frac{(-1)^{k+1}}{\prod_{\ell=0}^{k}w(\ell)}\sum_{i\geq 0}s^w(k+1,i)c(n-1+i)\right)$$

$$= \frac{(-1)^k}{\prod_{\ell=0}^{k-1}w(\ell)}\sum_{i\geq 0}\big(w(k)s^w(k,i)+s^w(k+1,i)\big)c(n-1+i)$$

$$= \frac{(-1)^k}{\prod_{\ell=0}^{k-1}w(\ell)}\sum_{i\geq 1}s^w(k,i-1)c(n-1+i) = \frac{(-1)^k}{\prod_{\ell=0}^{k-1}w(\ell)}\sum_{i\geq 0}s^w(k,i)c(n+i) = m(n,k).$$

We used the fact that $w(k)s^w(k,0)+s^w(k+1,0)=0$.

For $w(n)=n+1$ and $a(n)=\dfrac{1}{n+1}$ this reduces to the original Akiyama-Tanigawa algorithm for $(b(n))$ as shown in [11].

Now we give a list of some interesting formulas.

**6.2.** For $w(n)=n+1$, $a(n)=\dfrac{1}{n+1}$ we get from (2.16) and (3.9)

$$\sum_{j=0}^{n}S(n+1,j+1)(-1)^j\frac{j!}{j+1}=b(n). \tag{6.6}$$

Another proof of (6.6) can be found in [11], but I suspect that this result must be much older.

Formula (6.3) gives

$$\sum_{j=0}^{n}s(n+1,j+1)b(j)=(-1)^n\frac{n!}{(n+1)}. \tag{6.7}$$

**6.3.** For $w(n)=(n+1)^2$, $a(n)=n+1$ we know from (4.7) and (4.16) that $c(n)=(-1)^nG_{2n+2}$.

This gives

$$(-1)^{n-1}G_{2n}=\sum_{k=1}^{n}(-1)^{k-1}T(n,k)k\big((k-1)!\big)^2 \tag{6.8}$$

and

$$\sum_{k=1}^{n}(-1)^{n-k}t(n,k)G_{2k} = n!(n-1)!. \tag{6.9}$$

The Akiyama-Tanigawa algorithm applied to (6.8) gives another method for computing the Genocchi numbers. Choose $w(n) = (n+1)^2$ and $a(n) = n+1$ in Theorem 7.2.

Then the left upper part of the corresponding matrix $M = (m(i,j))$ is given by

$$\begin{pmatrix} 1 & 2 & 3 & 4 & 5 & 6 \\ -1 & -4 & -9 & -16 & -25 & -36 \\ 3 & 20 & 63 & 144 & 275 & 468 \\ -17 & -172 & -729 & -2096 & -4825 & -9612 \\ 155 & 2228 & 12303 & 43664 & 119675 & 276660 \\ -2073 & -40300 & -282249 & -1216176 & -3924625 & -10444428 \end{pmatrix}.$$

In the first column we get $c(n) = (-1)^n G_{2n+2}$.

**6.4.** From (4.28) we deduce the following identity (cf. [14], Exercise 5.8):

$$(-1)^{n-1}G_{2n+2} = \sum_{k=1}^{n}(-1)^{k-1}T(n,k)((k)!)^2. \tag{6.10}$$

The left upper part of the corresponding Akiyama-Tanigawa matrix is

$$\begin{pmatrix} 1 & 4 & 9 & 16 & 25 & 36 \\ -3 & -20 & -63 & -144 & -275 & -468 \\ 17 & 172 & 729 & 2096 & 4825 & 9612 \\ -155 & -2228 & -12303 & -43664 & -119675 & -276660 \\ 2073 & 40300 & 282249 & 1216176 & 3924625 & 10444428 \\ -38227 & -967796 & -8405343 & -43335184 & -162995075 & -495672372 \end{pmatrix}.$$

(6.3) gives the companion formula

$$\sum_{k=0}^{n}(-1)^{n-k}t(n,k)G_{2k+2} = (n!)^2. \tag{6.11}$$

For example for $n=3$ we have $4G_4 + 5G_6 + G_8 = 4+15+17 = 36 = (3!)^2$.

**6.5.** (4.21) and Lemma 4.3 give

$$\sum_{k=0}^{n}(-1)^{n-k}LS(n+1,k+1)((k+1)!)^2 = H_{2n+3}. \tag{6.12}$$

**6.6.** From (4.43) and (4.42) we deduce for $w(n) = (n+2)^2$ and $a(n) = n+2$

$$\sum_{k=0}^{n} (-1)^{n-k} S^w(n,k)(k+1)!(k+2)! = G_{2n+2} + G_{2n+4} \tag{6.13}$$

and

$$\sum_{k=0}^{n} (-1)^{n-k} s^w(n,k)\left(G_{2k+2} + G_{2k+4}\right) = (n+1)!(n+2)!. \tag{6.14}$$

**6.7.** For $w(n) = n+1$ and $a(n) = \dfrac{1}{n+1}$ we get from (4.47)

$$(2n+1)B_{2n} = \sum_{j=0}^{n} (-1)^j \frac{(j!)^2}{j+1} T(n+1, j+1). \tag{6.15}$$

Here the left upper part of the Akiyama-Tanigawa matrix begins with

$$\begin{pmatrix} 1 & \frac{1}{2} & \frac{1}{3} & \frac{1}{4} & \frac{1}{5} & \frac{1}{6} \\ \frac{1}{2} & \frac{2}{3} & \frac{3}{4} & \frac{4}{5} & \frac{5}{6} & \frac{6}{7} \\ -\frac{1}{6} & -\frac{1}{3} & -\frac{9}{20} & -\frac{8}{15} & -\frac{25}{42} & -\frac{9}{14} \\ \frac{1}{6} & \frac{7}{15} & \frac{3}{4} & \frac{104}{105} & \frac{25}{21} & \frac{19}{14} \\ -\frac{3}{10} & -\frac{17}{15} & -\frac{303}{140} & -\frac{16}{5} & -\frac{25}{6} & -\frac{353}{70} \\ \frac{5}{6} & \frac{433}{105} & \frac{261}{28} & \frac{232}{15} & \frac{460}{21} & \frac{4353}{154} \end{pmatrix}.$$

**6.8.** From (5.6) and (5.10) we deduce

$$\sum_{k=0}^{n} (-1)^{n-k} 4^{n-k} U(n,k)(2k+1)\left((2k-1)!!\right)^2 = T_{2n+1} \tag{6.16}$$

**6.9.** Finally Proposition 5.3 gives

$$\sum_{k=0}^{n} (-1)^k U(n,k) \frac{1}{(2k+1)4^k} \left((2k-1)!!\right)^2 = B_{2n}. \tag{6.17}$$